\documentclass[11pt,draft]{amsart}
\usepackage{amssymb,amscd,amsgen,amsxtra,amsthm}

\newtheorem{thm}{Theorem}[section]
\newtheorem{lem}[thm]{Lemma}
\newtheorem{cor}[thm]{Corollary}
\newtheorem{prop}[thm]{Proposition}

\newtheorem{defn}[thm]{Definition}
\theoremstyle{definition}
\newtheorem{rem}[thm]{Remark}
\newtheorem{expl}[thm]{Example}

\advance \topmargin by -\headheight
\advance \topmargin by -\headsep
     
\evensidemargin \oddsidemargin
\newcommand{\N}{\mathbf{N}}
\newcommand{\cs}{\mathcal{S}}
\newcommand{\cn}{\mathcal{N}}
\newcommand{\cf}{\mathcal{F}}

\begin{document}

\title{A dichotomy on Schreier sets}
\author{Robert Judd}
\address{Department of Mathematics\\
  University of Texas at Austin\\
  Austin, TX 78712-1082 \\
  U.S.A. }
\email{rjudd@math.utexas.edu}

\begin{abstract}
  We show that the Schreier sets $\mathcal{S}_{\alpha}\ 
  (\alpha<\omega_1)$ satisfy the following dichotomy property.  For
  every hereditary collection $\cf$ of finite subsets of $\N$, either
  there exists infinite $M=(m_i)_1^{\infty}\subseteq\N$ such that
  $\cs_{\alpha}(M)=\{\{m_i:i\in E\}:E\in\cs_{\alpha}\}\subseteq\cf$,
  or there exist infinite $M=(m_i)_1^{\infty},N\subseteq\N$ such that
  $\cf[N](M)=\{\{m_i:i\in F\}:F\in\cf \mbox{ and } F\subset
  N\}\subseteq\cs_{\alpha}$.
\end{abstract}

\maketitle

\section{Introduction}
\label{sec:intro}

Collections of finite subsets of the natural numbers have become
important in Banach space theory.  The Schreier sets $\cs_{\alpha}$,
defined below for each countable ordinal $\alpha$, are the most common
among these sets.  The first Schreier set, $\cs_1$, is fundamental to
the construction of the original Tsirelson space, see~\cite{t}
and~\cite{fj}, while the more general Schreier sets are used to
construct the Schreier spaces, which may be found
in~\cite{sh},~\cite{aa} and~\cite{ao}, and the exciting new collection
of Tsirelson type spaces developed by Argyros and Deliyanni~\cite{ad}.

The Banach spaces mentioned above may be constructed with collections
of finite subsets of the natural numbers other than the Schreier sets.
However, the Schreier sets are in some sense universal for these
alternate collections.  For example a result of Odell, Tomczak and
Wagner~\cite{otw} shows that for pointwise closed collections $\cf$ of
finite subsets of $\N$ there exists a subsequence $M$ of $\N$ such
that $\cf(N)$ is a subset of one of the Schreier sets.  (The notation
$\cf(N)$ is described below.)

We show roughly that if we fix a Schreier set $\cs_{\alpha}$, then
herediary collections $\cf$ of finite subsets of $\N$ satisfy: either
$\cf$ is of sufficient complexity to contain the Schreier set or the
sets in $\cf$ lying in some subsequence must be contained in the
Schreier set.  The precise statement is a bit more complicated.  One
must allow for a wide range of collections of finite subsets.  For
example the first Schreier set, $\cs_1$, consists of all finite
subsets of $\N$ whose smallest element is at least as large as the
size of the set.  This condition is called an admissibility condition.
There are many such conditions.  A different admissibility condition
would be to consider collections of finite subsets such that the
square of the smallest element in each set is at least as large as the
number of elements in the set.  We state the dichotomy theorem here,
deferring the notation until Section~\ref{sec:prelims}.

\begin{thm} \label{mainthm}
  For each $\alpha<\omega_1$, for every hereditary collection
  $\cf\subseteq[\N]^{<\omega}$ and for all $\bar M\in[\N]$ either
  there exists $M\in[\bar M]$ such that $\cs_{\alpha}(M)\subseteq\cf$
  or there exist $M\in[\bar M],\ N\in[\N]$ such that
  $\cf[N](M)\subseteq\cs_{\alpha}$.
\end{thm}

In the next section we define the Schreier classes $\cs_{\alpha}\
(\alpha<\omega_1)$ along with other notions concerning collections of
finite subsets of $\N$.  We also introduce Schreier games; these are a
method of choosing finite subsets of $\N$ in such a way that the
resulting set is in one of the collections $\cs_{\alpha}$.  The
combinatorial framework for proving Theorem~\ref{mainthm} is presented
in Section~\ref{sec:dichotp} as the dichotomy property.  We devote
Section~\ref{sec:app} to an alternative proof of a result of Argyros,
Mercourakis and Tsarpalias~\cite{amt}, using Theorem~\ref{mainthm}.

\section{Preliminaries}
\label{sec:prelims}

We use various subsets, and collections of subsets, of the natural
numbers throughout this paper; for future reference we define all the
notation for these sets at the beginning of this section.  In general
$L,\ M$ and $N$ will be infinite subsets of $\N$, while $E$ and $F$
will be finite subsets, and $\cf$ and $\mathcal{G}$ will be
collections of finite subsets of $\N$.  We consider every subset of
$\N$, whether finite or infinite, to be an increasing sequence.  Thus
if $N\subseteq\N$, then $N=(n_i)_{i=1}^{\infty}$ where
$n_1<n_2<\dots$ and if $E$ is a finite subset, then
$E=\{e_1,\dots,e_k\}$ where $e_1<\dots<e_k$.

When $N$ is an infinite subset of $\N$ we let $[N]$ be the set of
infinite subsets of $N$ and we let $[N]^{<\omega}$ represent the set
of finite subsets of $N$.  Let $E,F\subseteq{\N}$ and $n\geq1$.  We
write $E<F$ if either set is empty or if $\max E<\min F$, $n<E$ if
$\{n\}<E$ and $n\leq E$ if $n\leq\min E$.

Let $\cf$ be a collection of finite subsets of $\N$.  We next define
three properties which $\cf$ may have: hereditary, spreading and
closed.  For $\cf$ to be \emph{hereditary} requires that whenever
$E\subset F$ and $F\in\cf$ then $E\in\cf$.  We say $\cf$ is
\emph{spreading} if whenever $F=\{m_1,\dots,m_k\}\in\cf$ and
$n_1<\dots<n_k$ satisfies: $m_i\leq n_i$ for $i\leq k$, then
$\{n_1,\dots,n_k\}\in\cf$.  The set $2^{\N}$ of all subsets of $\N$
is a topological space under the topology of pointwise convergence;
thus $\cf$ is \emph{(pointwise) closed} if it is closed in $2^{\N}$.
We collect the first and third properties together and say that $\cf$
is \emph{adequate} if it is both closed and hereditary.

Finally we need some notation to talk about what happens when we
restrict a collection of finite subsets of $\N$ to an infinite subset
of $\N$.  Let $N=(n_i)\in[\N]$ be an infinite sequence and let
$\cf\subseteq[\N]^{<\omega}$ be a collection of finite subsets.  We
write the subset of $\cf$ consisting only of those elements which are
also subsets of $N$ as $\cf[N]$.  Thus
\[ \cf[N] = \{ F\in\cf : F\subseteq N \} = 
\cf \cap [N]^{<\omega}\ . \]
We also want to put $\cf$ into the sequence $N$.  In other words if
$F\in\cf$ and we define $n_F=\{ n_i : i\in F \}$, then $\cf(N)$ is the
collection of all such sets, ie.\ $\cf(N)=\{ n_F : F\in \cf \}$.
(Note that $n_F$ is a finite subset of $N$.)

\begin{defn}\label{def:schreier}
  The Schreier sets, $\mathbf{\cs_{\alpha}}$ \emph{\cite{aa}}   
\end{defn}
\noindent\nobreak
The Schreier sets, $\cs_{\alpha}$ for each $\alpha<\omega_1$, are
defined inductively as follows: let $\cs_0=\{\{n\} :
n\geq1\}\cup\{\emptyset\}$ and $\cs_1=\{F\subset {\mathbf N} : |F|\leq
F\}$.  (Note that this definition allows for $\emptyset\in\cs_1$.)  If
$\cs_{\alpha}$ has been defined let
\[\cs_{\alpha+1}=\{ \cup_1^kF_i : k\leq
F_1<\dots<F_k,\ F_i\in \cs_{\alpha}\ (i=1,\dots,k),\ k\in\N \}\ .\]
If $\alpha$ is a limit ordinal with $\cs_{\beta}$ defined for each
$\beta<\alpha$, choose and fix an increasing sequence of ordinals
$(\alpha_n)$ with $\alpha=\sup_n\alpha_n$ and let
\[\cs_{\alpha}=
\cup_{n=1}^{\infty}\{ F\in\cs_{\alpha_n} : n\leq F\}\ .\]
Note that each $\cs_{\alpha}$ is hereditary, spreading and closed.
For $r\geq1$ and $\alpha_1,\dots,\alpha_r<\omega_1$ let
\[(\cs_{\alpha_1},\dots,\cs_{\alpha_r})=\{F=\cup_1^rF_i :
F_i\in\cs_{\alpha_i}\ (i\leq r)\mbox{ and }F_1<\dots<F_r\}\ .\]

\begin{defn}Schreier games\end{defn}
\noindent\nobreak
We define a game for two players on $\N$, called an
\emph{$(\alpha_1,\dots,\alpha_r)$-Schreier game}, for each $r$-tuple
of ordinals with $0\leq\alpha_1\leq\dots\leq\alpha_r<\omega_1$.  If
$r=1$, then we drop the parentheses and simply call it an
$\alpha$-Schreier game.  The two players are $\cn$ who chooses numbers
and $\cs$ who chooses non-empty sets.  Roughly, $\cn$ will pick a
finite sequence of numbers and $\cs$ will pick a finite sequence of
finite subsets of $\N$, $E_1<\dots<E_k$.  The number of choices made
and the order of the plays will depend upon the particular
$(\alpha_1,\dots,\alpha_r)$-Schreier game being played, and may also
depend upon previous plays.

We first describe the $\alpha$-Schreier game for $\alpha<\omega_1$.
In the \emph{0-Schreier game} $\cs$ chooses $\{n\}$ for some $n\geq1$.
In the 1-Schreier game $\cn$ picks $l\geq1$ and $\cs$ chooses $E\in
[\N]^{<\omega}$ such that $|E|\geq l$.  Suppose we have already
described the $\alpha$-Schreier game for $\alpha<\omega_1$.  The
\emph{$(\alpha+1)$-Schreier game} starts with $\cn$ picking $l\geq1$
and then the two players play the $\alpha$-Schreier game $l$ times,
with the additional condition that if $E$ is the last set $\cs$ chose
in the $i^{\,\mbox{\scriptsize th}}\ \alpha$-Schreier game and $F$ is
the first set $\cs$ chose in the $(i+1)^{\mbox{\scriptsize th}}\
\alpha$-Schreier game, then $E<F$.  For $\alpha$ a limit ordinal
suppose we have already described the $\gamma$-Schreier game for each
$\gamma<\alpha$ and let $\alpha_n\nearrow\alpha$ be the sequence used
to define $\cs_{\alpha}$.  The $\alpha$-Schreier game starts with
$\cn$ picking $l\geq1$ and then the two players play the
$\alpha_l$-Schreier game.

If $\alpha_1\leq\dots\leq\alpha_r<\omega_1$, then an
\emph{$(\alpha_1,\dots,\alpha_r)$-Schreier game} is simply an
$\alpha_1$-Schreier game followed by an $\alpha_2$-Schreier game, and
so on, finishing with an $\alpha_r$-Schreier game.  The only other
condition is that if $E$ is the last set $\cs$ chose in the
$\alpha_i$-Schreier game and $F$ is the first set $\cs$ chose in the
$\alpha_{i+1}$-Schreier game, then $E<F$.  In the sequel, by an
\emph{$(\alpha_1,\dots,\alpha_r)$-game} we shall mean an
$(\alpha_1,\dots,\alpha_r)$-Schreier game.

As an example, consider the 2-Schreier game.  $\cn$ chooses $l\geq1$
and then they play the $(1\mbox{,$\,\stackrel{l}{\dots}\,$,}1)$-game.
This starts with $\cn$ choosing $k_1\geq1$ and then $\cs$ chooses
$E_1$ with $|E_1|\geq k_1$.  Then $\cn$ chooses $k_2$ and $\cs$
chooses $E_2$ with $|E_2|\geq k_2$ and $E_2>E_1$.  This continues
until $\cn$ has chosen $k_l$ and $\cs$ has chosen $E_l$ with
$|E_l|\geq k_l$ and $E_l>E_{l-1}$.  The set resulting from this game
is $E=\cup_1^lE_i$.  In general if $\cn$ and $\cs$ play a Schreier
game, and $(E_i)^k_1$ is the sequence of sets which $\cs$ chose in the
game, with $E_i<E_{i+1}\ (1\leq i<k)$, then \emph{the set $E$
resulting from the game} is defined as $E=\cup_1^kE_i$.

A \emph{bound} $(\alpha_1,\dots,\alpha_r)$-Schreier game is one where
at each stage $\cn$ is restricted to exactly one choice of number to
pick.  If $\cn$ and $\cs$ play a bound game and $E$ is the set
resulting from this game, then we say \emph{$\cs$ chose $E$ as small
  as possible} if at each stage, when $\cs$ had to choose a set $E_i$
of size at least $l_i$, then $\cs$ always chose $E_i$ of size equal to
$l_i$.

We say that $\cn$ has a \emph{winning strategy for the
  $(\alpha_1,\dots,\alpha_r)$-Schreier game on
  $\cf\subseteq[\N]^{<\omega}$} if $\cn$ can choose integers so that,
  whatever sets $\cs$ picks, the set $E$ resulting from the game does
  not belong to $\cf$.  Notice that if $\cn$ has a winning strategy
  for the $(\alpha_1,\dots,\alpha_r)$-Schreier game on
  $\cf\subseteq[\N]^{<\omega}$, then $\cn$ also has a winning strategy
  for the $(\alpha_1,\dots,\alpha_r)$-Schreier game on $\cf[M]$ for
  any $M\in[\N]$ since $\cf[M]\subseteq\cf$.

As an example of a winning strategy for $\cn$ we shall consider the
$(1,1)$-game on $\cs_1$.  In this game $\cn$ chooses $l\geq1$, next
$\cs$ chooses $E\in[\N]^{<\omega}$ with $|E|\geq l$, then $\cn$
chooses $m\geq1$ and finally $\cs$ chooses $F>E$ with $|F|\geq m$.  A
winning strategy for $\cn$ in this game would be to choose $l=1$ and
$m=\min E$ (which $\cn$ may do since $\cs$ chooses $E$ before $\cn$
chooses $m$).  Now, if $A=E\cup F,$ then $|A|\geq l+m=1+m>\min A$,
while if $A$ were in $\cs_1$, then we would have $|A|\leq\min A$.
Thus $A\not\in\cs_1$ which is what $\cn$ was trying to achieve.

\section{The Dichotomy property}
\label{sec:dichotp}

\begin{defn}The Dichotomy property, (D)\end{defn}
\noindent
An $r$-tuple of ordinals, $(\alpha_1,\dots,\alpha_r)$ with
$0\leq\alpha_1\leq\dots\leq\alpha_r<\omega_1$, has the Dichotomy
property (D) if for each hereditary collection
$\cf\subseteq[\N]^{<\omega}$ and every $\bar N\in[\N]$, either there
exists $M\in[\bar N]$ such that
$(\cs_{\alpha_1},\dots,\cs_{\alpha_r})(M)\subseteq\cf$, or there
exists $M\in[\bar N]$ such that $\cn$ has a winning strategy for the
$(\alpha_1,\dots,\alpha_r)$-Schreier game on $\cf[M]$.

This section is devoted to proving that every increasing $r$-tuple of
countable ordinals has the Dichotomy property.

\begin{prop} \label{p:dichot}
  The $r$-tuple $(\alpha_1,\dots,\alpha_r)$ has the Dichotomy property
  (D) for each $r\geq1$ and every $r$-tuple of ordinals with
  $0\leq\alpha_1\leq\dots\leq\alpha_r<\omega_1$.
\end{prop}

We prove this inductively in several stages using a technique
developed by Kiriakouli and Negrepontis~\cite{kn}.  The method
consists of a double induction.  To prove that every $r$-tuple of
ordinals, $(\alpha_1,\dots,\alpha_r)$, has a certain property (P) one
first shows that if $(\alpha_1,\dots,\alpha_r)$ has (P), then so does
$(\alpha,\alpha_1,\dots,\alpha_r)$.  Next one demonstrates that if
$(\alpha\mbox{,$\,\stackrel{k}{\dots}\,$,}
\alpha,\alpha_1,\dots,\alpha_r)$ has the property for every $k\geq1$,
then so does $(\alpha+1,\alpha_1,\dots,\alpha_r)$.  The rest of the
proof usually follows easily from these two results.  In our case the
key to proving Proposition~\ref{p:dichot} is the following lemma:

\begin{lem} \label{l:first}
  Let $r\geq1$ and let $0\leq\alpha_1\leq\dots\leq\alpha_r<\omega_1$.
  If $(\alpha_1,\dots,\alpha_r)$ has the Dichotomy property (D), then
  so does $(0,\alpha_1,\dots,\alpha_r)$.
\end{lem}
\begin{proof}
  Let $\cf\subseteq[\N]^{<\omega}$ be hereditary and $\bar N=(\bar
  n_i)_{i=1}^{\infty}\in[\N]$, then we seek $L\in[\bar N]$ such that
  either $(\cs_0,\cs_{\alpha_1},\dots,\cs_{\alpha_r})(L)\subseteq\cf$
  or $\cn$ has a winning strategy for the
  $(0,\alpha_1,\dots,\alpha_r)$-game on $\cf[L]$.  We cannot find $L$
  all at once; instead we must choose it bit by bit.  We construct
  sequences $M_l=(m^l_i)_{i=1}^{\infty}$ with $\bar N=M_0\supseteq
  M_1\supseteq M_2\supseteq\dotsb$ such that either $\{m_l^l\}\cup
  F\in\cf$ for each $F\subseteq\cs(M_l)$ with $F>m^l_l$, or else $\cn$
  has a winning strategy in the $(0,\alpha_1,\dots,\alpha_r)$-game on
  $\cf[M_l]$ provided the first choice of $\cs$ is $m_l^l$.  We may
  then choose $L$ as a diagonal subsequence of these sequences $M_l$.

  We begin by defining
  \[\cf_1=\{ F : \{\bar n_1\}\cup F\in\cf[M_0] \}\ .\]
  Since $(\alpha_1,\dots,\alpha_r)$ has (D), it follows that there
  exists $\bar M_1=(\bar m^1_i)_{i=1}^{\infty}\in[\bar N]$ such that
  either $\cn$ has a winning strategy for the
  $(\alpha_1,\dots,\alpha_r)$-Schreier game on $\cf_1[\bar M_1]$ or
  $(\cs_{\alpha_1},\dots,\cs_{\alpha_r})(\bar M_1)$ is a subset of
  $\cf_1$.  Let $M_1=(\bar n_1,\bar m^1_2,\bar m^1_3,\dots)$ be the
  sequence $\bar M_1$ with its first element replaced by $\bar n_1$.
  Now, either $\{\bar n_1\}\cup E\in\cf$ for each
  $E\in(\cs_{\alpha_1},\dots,\cs_{\alpha_r})(M_1)$ with $E>\bar n_1$
  or $\{\bar n_1\}\cup F\not\in\cf$ for every set $F\subseteq
  M_1\setminus\{\bar n_1\}$ resulting from $\cn$ playing a winning
  strategy for the $(\alpha_1,\dots,\alpha_r)$-game on $\cf_1[M_1]$.
  This last follows since if $F\in\cf_1[M_1]$ and $F>\bar n_1$, then
  $F\in\cf_1[\bar M_1]$.
  
  Suppose we have chosen sequences $\bar N\supseteq M_1\supseteq
  M_2\supseteq\dots\supseteq M_{l-1}$ with the properties:
  \begin{itemize}
  \item If $M_i=(m^i_j)_{j=1}^{\infty}$ for $1\leq i<l$, then
    $m_j^{i-1}=m_j^i$ whenever $1\leq j\leq i$ and $1<i<l$.
  \item For each $i=1,\dots,l-1$ either $\{m_i^i\}\cup F\in\cf$ for
    all $F\in(\cs_{\alpha_1},\dots,\cs_{\alpha_r})(M_i)$ with
    $F>m_i^i$, or else $\{m_i^i\}\cup F\not\in\cf$ for any
    $F\subseteq(m_j^i)_{j=i+1}^{\infty}$ resulting from $\cn$ playing
    a winning strategy in the $(\alpha_1,\dots,\alpha_r)$-game on
    $\cf_i[M_i]$, where $\cf_i=\{F:\{m_i^{i-1}\}\cup F
    \in\cf[M_{i-1}]\}$.
  \end{itemize}
  
  To construct the next sequence $M_l$ we define
  \[ \cf_l=\{F:\{m^{l-1}_{l}\}\cup F\in\cf[M_{l-1}]\}\ .\]  
  Since $(\alpha_1,\dots,\alpha_r)$ has (D), it follows that there
  exists $\bar M_{l}=(\bar m^{l}_i)_i\in[M_{l-1}]$ such that either
  $\cn$ has a winning strategy for the
  $(\alpha_1,\dots,\alpha_r)$-Schreier game on $\cf_{l}[\bar M_{l}]$
  or $(\cs_{\alpha_1},\dots,\cs_{\alpha_r})(\bar
  M_{l})\subseteq\cf_{l}$.  Let
  \[M_{l}=(m^{l-1}_1,\dots,m^{l-1}_{l},\bar m^{l}_{l+1},\bar 
  m^{l}_{l+2},\dots)\] 
  be the sequence $\bar M_l$ with the first $l$
  elements replaced by the first $l$ elements of $M_{l-1}$.  As with
  $M_1$ and $\cf_1$, either $\{m^{l}_{l}\}\cup F\in\cf$ for each
  $F\in(\cs_{\alpha_1},\dots,\cs_{\alpha_r})(M_{l})$ with
  $F>m^{l}_{l}$ or $\{m^{l}_{l}\}\cup F\not\in\cf$ for every
  $F\subseteq(m^l_i)_{i>l}$ resulting from $\cn$ playing a winning
  strategy for the $(\alpha_1,\dots,\alpha_r)$-game on
  $\cf_{l}[M_{l}]$.
  
  We repeat this process for each $l\geq1$.  Let
  $M=(m_k)_{k=1}^{\infty}$ be the sequence defined by $m_k=m^k_k$ for
  each $k\geq1$.  Then for each $l\geq1$ either $\{m_l\}\cup F\in\cf$
  for all $F\in(\cs_{\alpha_1},\dots,\cs_{\alpha_r})(M)$ with $F>m_l$
  or $\{m_l\}\cup F\not\in\cf$ for each $F\subseteq(m_k)_{k>l}$
  resulting from $\cn$ playing a winning strategy for the
  $(\alpha_1,\dots,\alpha_r)$-game on $\cf_l[M]$.  This induces a
  coloring on $\N$; in the first case we color $l\in\N$ red, and in
  the second, blue.
  
  Now, either there exists an infinite subsequence $J\in[\N]$ such
  that every $j\in J$ is colored red, in which case let
  $L=(m_{j})_{j\in J}$, or there exists $k\geq1$ such that $l$ is
  colored blue for all $l\geq k$, and let $L=(m_l)_{l\geq k}$.  In the
  first case it is clear that
  $(\cs_0,\cs_{\alpha_1},\dots,\cs_{\alpha_r})(L)\subseteq\cf$.  In
  the second case if $\cs$ picks $\{n\}$ with $n\not\in L$, then the
  resulting set cannot be in $\mathcal{F}[L]$.  Otherwise $\cs$ picks
  $\{m_l\}$ for some $l\geq k$ and then $\cn$ has a winning strategy
  for the $(\alpha_1,\dots,\alpha_r)$-Schreier game on
  $\mathcal{F}_l[(m_i)_{i>l}]$.  In either situation we see that $\cn$
  has a winning strategy for the
  $(0,\alpha_1,\dots,\alpha_r)$-Schreier game on $\mathcal{F}[L]$ as
  required.
\end{proof}

\begin{lem}\label{l:aaa=>a+1}
  If $0\leq\alpha<\alpha_1\leq\dots\leq\alpha_r$, for some $r\geq1$,
  and the $(k+r)$-tuple $(\alpha\mbox{,$\,\stackrel{k}{\dots}\,$,}
  \alpha,\alpha_1,\dots,\alpha_r)$ has property (D) for every
  $k\geq1$, then $(\alpha+1,\alpha_1,\dots,\alpha_r)$ has property
  (D).
\end{lem}
\begin{proof}
  Let $\cf\subseteq[\N]^{<\omega}$ be hereditary, let $\bar N\in[\N]$
  and find sequences $\bar N\supseteq L_1\supseteq L_2\supseteq\dotsb$
  such that for each $k$ either
  $(\cs_{\alpha}\mbox{,$\,\stackrel{k}{\dots}\,$,}
  \cs_{\alpha},\cs_{\alpha_1},\dots,\cs_{\alpha_r})(L_k)
  \subseteq\cf$ or $\cn$ has a winning strategy in the
  $(\alpha\mbox{,$\,\stackrel{k}{\dots}\,$,}
  \alpha,\alpha_1,\dots,\alpha_r)$-Schreier game on $\cf[L_k]$.  In
  this last case $\cn$ has a winning strategy in the
  $(\alpha+1,\alpha_1,\dots,\alpha_r)$-Schreier game on $\cf[L_k]$
  given by: for the $(\alpha+1)$-game $\cn$ picks $k$ and then plays a
  winning strategy in the $(\alpha\mbox{,$\,\stackrel{k}{\dots}\,$,}
  \alpha,\alpha_1,\dots,\alpha_r)$-Schreier game on $\cf[L_k]$.
  Otherwise we set $L=(l^k_k)$ and then we obtain
  $(\cs_{\alpha+1},\cs_{\alpha_1},\dots,\cs_{\alpha_r})(L)
  \subseteq\cf$.
\end{proof}

\begin{lem}\label{l:b<a=>a}
  If $\alpha$ is a limit ordinal with
  $\alpha\leq\alpha_1\leq\dots\leq\alpha_r$, for some $r\geq1$, and
  $(\beta,\alpha_1,\dots,\alpha_r)$ has property (D) for every
  $\beta<\alpha$, then $(\alpha,\alpha_1,\dots,\alpha_r)$ has property
  (D).
\end{lem}
\begin{proof} 
  As before, let $\cf\subseteq[\N]^{<\omega}$ be hereditary, let $\bar
  N\in[\N]$ and find sequences $\bar N\supseteq L_1\supseteq
  L_2\supseteq\dotsb$ such that for every $k$ either
  $(\cs_{\beta_k},\cs_{\alpha_1},\dots,\cs_{\alpha_r})(L_k)
  \subseteq\cf$ or $\cn$ has a winning strategy in the
  $(\beta_k,\alpha_1,\dots,\alpha_r)$-Schreier game on $\cf[L_k]$
  where $\beta_k\nearrow\alpha$ is the sequence used to define
  $\cs_{\alpha}$.  In the second case $\cn$ has a winning strategy in
  the $(\alpha,\alpha_1,\dots,\alpha_r)$-game on $\cf[L_k]$ by
  choosing $k$ and playing a winning strategy in the
  $(\beta_k,\alpha_1,\dots,\alpha_r)$-Schreier game on $\cf[L_k]$.
  Otherwise let $L=(l^k_k)$, then
  $(\cs_{\alpha},\cs_{\alpha_1},\dots,\cs_{\alpha_r})(L)
  \subseteq\cf$.
\end{proof}

\begin{lem}\label{l:a1-ar=>b1-bsa1-ar}
  If $0\leq\alpha_1\leq\dots\leq\alpha_r$, for some $r\geq1$, and
  $(\alpha_1,\dots,\alpha_r)$ has property (D) then so too does
  $(\beta_1,\dots,\beta_s,\alpha_1,\dots,\alpha_r)$ for all
  $\beta_1\leq\dots\leq\beta_s\leq\alpha_1$, and each $ s\geq1$.
\end{lem}
\begin{proof}
  We prove this by induction on $\beta_s$ for arbitrary $s$ and
  $\alpha_1\geq\beta_s$. When $\beta_s=0$ the result follows by
  iterating Lemma~\ref{l:first} $s$-times.
  
  Suppose that we have proven the result for $\beta=\beta_s$, ie.\ we
  have shown that for every $\alpha_r\geq\dots\geq\alpha_1\geq\beta$,
  if $(\alpha_1,\dots,\alpha_r)$ has property (D), then so too does
  $(\gamma_1,\dots,\gamma_k,\beta,\alpha_1,\dots,\alpha_r)$ for all
  $\gamma_1\leq\dots\leq\gamma_k\leq\beta$.  Clearly we may take each
  $\gamma_i=\beta$ and so in particular we have proven that for every
  $\alpha_r\geq\dots\geq\alpha_1>\beta$, if
  $(\alpha_1,\dots,\alpha_r)$ has property (D), then
  $(\beta\mbox{,$\,\stackrel{k}{\dots}\,$,}
  \beta,\alpha_1,\dots,\alpha_r)$ also has property (D) for each
  $k\geq1$.  Hence by Lemma~\ref{l:aaa=>a+1} so does
  $(\beta+1,\alpha_1,\dots,\alpha_r)$.  Iterating this argument, we
  obtain that if $\alpha_r\geq\dots\geq\alpha_1>\beta$ and
  $(\alpha_1,\dots,\alpha_r)$ has property (D), then
  $(\beta+1\mbox{,$\,\stackrel{l}{\dots}\,$,}
  \beta+1,\alpha_1,\dots,\alpha_r)$ also has property (D) for all
  $l\geq1$, and hence $(\gamma_1,\dots,\gamma_k,
  \beta+1\mbox{,$\,\stackrel{l}{\dots}\,$,}
  \beta+1,\alpha_1,\dots,\alpha_r)$ does too, for all
  $\gamma_1\leq\dots\leq\gamma_k\leq\beta$ by the result for $\beta$.
  In other words, for every
  $\alpha_r\geq\dots\geq\alpha_1\geq\beta+1$, if
  $(\alpha_1,\dots,\alpha_r)$ has property (D), then so too does
  $(\beta_1,\dots,\beta_s,\alpha_1,\dots,\alpha_r)$ for all
  $\beta_1\leq\dots\leq\beta_s\leq\beta+1$, and any $ s\geq1$ as
  required.
  
  If $\beta_s$ is a limit ordinal and we have proven the result for
  each $\beta<\beta_s$, then $(\beta,\alpha_1,\dots,\alpha_r)$ has
  property (D) for each $\beta<\beta_s$.  Thus, by
  Lemma~\ref{l:b<a=>a}, the $(r+1)$-tuple
  $(\beta_s,\alpha_1,\dots,\alpha_r)$ also has property (D).  Now, as
  in the successor case, we have
  that$(\beta_1,\dots,\beta_s,\alpha_1,\dots,\alpha_r)$ has property
  (D) as required.
\end{proof}

\begin{proof}[Proof of Proposition~\ref{p:dichot}]
  We prove by induction that $(\alpha)$ has property (D) for each
  $\alpha<\omega_1$, and then the result follows from
  Lemma~\ref{l:a1-ar=>b1-bsa1-ar}.  Let $\alpha=0$, let
  $\cf\subseteq[\N]^{<\omega}$ be hereditary and let $\bar N\in[\N]$.
  Let $L=\{n\in\bar N : \{n\}\in\cf\}$; if $L$ is infinite let $M=L$,
  and then $\cs_0(M)\subseteq\cf$.  Otherwise let $M=\bar N\setminus
  L$, then $\cn$ has a winning strategy for the 0-Schreier game on
  $\cf[M]$ since this set is empty.  This completes the proof for
  $\alpha=0$.
  
  If $(\alpha)$ has (D) then by Lemma~\ref{l:a1-ar=>b1-bsa1-ar} so
  does $(\alpha\mbox{,$\,\stackrel{k}{\dots}\,$,}\alpha)$ for each
  $k\geq1$.  Thus $(\alpha+1)$ has (D) by Lemma~\ref{l:aaa=>a+1}.  If
  $\alpha$ is a limit ordinal and $(\beta)$ has (D) for each
  $\beta<\alpha$ then $(\alpha)$ has (D) by Lemma~\ref{l:b<a=>a}.
  This completes the proof.
\end{proof}

\section{The Main Result}

In this section we prove Theorem~\ref{mainthm}.  Actually we prove a
somewhat stronger statement:
\begin{thm}\label{bigthm}
  For all $r\geq1$ and each $r$-tuple of countable ordinals
  $0\leq\alpha_1\leq\dots\leq\alpha_r<\omega_1$, if $\cf$ is a
  hereditary collection, $\cf\subseteq[\N]^{<\omega}$ and $\bar
  N\in[\N]$, then either there exists $M\in[\bar N]$ such that
  $(\cs_{\alpha_1},\dots,\cs_{\alpha_r})(M)\subseteq\cf$, or there
  exist $M\in[\bar N],\ N\in[\N]$ such that
  $\cf[M](N)\subseteq(\cs_{\alpha_1},\dots,\cs_{\alpha_r})$.
\end{thm}

\begin{prop} \label{P:schreier}
  If $\cs$ and $\cn$ play a bound $(\alpha_1,\dots,\alpha_r)$-game,
  then there exists $N=(n_i)\in[\N]$ such that if $E$ is any result of
  this bound game where $\cs$ chooses $E$ as small as possible, then
  $n_E=\{ n_i : i\in E \}\in(\cs_{\alpha_1},\dots,\cs_{\alpha_r})$.
\end{prop}

Before we give the proof of this proposition we recall the notion of
spreading.  A collection $\cf\subseteq[\N]^{<\omega}$ is
\emph{spreading} if it has the property that if $G=\{ g_1,\dots,g_n
\}\in\cf$ and $H=\{ h_1,\dots,h_n \}$ satisfies: $g_j\leq h_j\ 
(j=1,\dots,n)$, then also $H\in\cf$.  In this case we say that $H$ is
a \emph{spreading} of $G$.  Moreover it is easy to see that if $\cf$
is spreading and $M=(m_i),\ N=(n_i)\in[\N]$ satisfy $m_i\leq n_i$ for
all $i$, then $m_E\in\cf$ implies that $n_E\in\cf$.

\begin{proof}
  We first prove the result for $r=1$ by induction on $\alpha$. This
  is then easy to generalize.  In order to find the sequence $N=(n_t)$
  we construct an increasing function $f:\N\rightarrow\N$ and let
  $n_t=f(t)$.
  
  \textit{Case 1, $\alpha=0$.} This is clearly true, just by setting
  $f(t)=t$.
  
  \textit{Case 2, $\alpha\Rightarrow\alpha+1$.} We assume that for any
  bound $\alpha$-game there exists a function $f$, as above.  Now, an
  $(\alpha+1)$-game consists of $\cn$ choosing $k$ and then the two
  players play a bound
  $(\alpha\mbox{,$\,\stackrel{k}{\dots}\,$,}\alpha)$-game.  Since the
  $(\alpha+1)$-game is bound there is only one choice of $k$ which
  $\cn$ may make.  For each of the bound $\alpha$-games which make up
  the $(\alpha+1)$-game we shall choose below a function $f^i$ such
  that $f^i(E_i)=\{ f^i(t) : t\in E_i \}\in\cs_{\alpha}$, for any set
  $E_i$ resulting from the $i^{\mbox{\scriptsize th}}\ \alpha$-game.
  We then let $f(t)=k+\sum_{i=1}^kf^i(t)$.
  
  The first $\alpha$-game is already fixed, so we may choose $f^1$
  using the hypothesis.  However the $i^{\mbox{\scriptsize th}}\
  \alpha$-game, while bound, depends on which sets were picked in the
  first $(i-1)$ games, so we cannot just pick $f^i$ straight from the
  hypothesis---instead we have to cover all possible bound
  $\alpha$-games which may be played.  Fortunately, for any fixed
  integer $t$ only finitely many bound
  $(\alpha\mbox{,$\,\stackrel{i-1}{\dots}\,$,}\alpha)$-games can be
  played which finish before $t$; let this number be $s$.  Thus there
  are $s$ possible bound $\alpha$-games we could be playing.  If
  $s\ne0$, then let the functions from these be $f^i_1,\dots,f^i_s$,
  and let $f^i(t)=\sum_{j=1}^sf^i_j(t)$.  Otherwise let $f^i(t)=t$.
  
  We must now show that the function $f$ given by
  $f(t)=k+\sum_{i=1}^kf^i(t)$ is the function we seek for the bound
  $(\alpha+1)$-game.  Let $E=\cup_1^kE_i$ be the result of the bound
  $(\alpha\mbox{,$\,\stackrel{k}{\dots}\,$,}\alpha)$-game where $\cs$
  chooses $E$ as small as possible and where $E_i$ is the result of
  the $i^{\mbox{\scriptsize th}}\ \alpha$-game.  We show that
  $f^i(E_i)\in\cs_{\alpha}$ for $i\leq k$.  We already know that this
  works for $i=1$ by the hypothesis.  Then for $1<i\leq k$, once we
  have chosen $E_1<\dots<E_{i-1}$ we will have fixed the bound
  $\alpha$-game we are playing when choosing $E_i$.  Let the function
  for this game be $f'$, from the induction hypothesis, then
  $f'(E_i)=\{ f'(t) : t\in E_i \}\in\cs_{\alpha}$.  But by the
  construction of $f^i$ we know that $f^i(t)\geq f'(t)$ for each $t$
  in $E_i$.  To obtain that $f^i(E_i)\in\cs_{\alpha}$ recall that the
  collection $\cs_{\alpha}$ is spreading and clearly $f^i(E_i)$ is a
  spreading of $f'(E_i)$, hence $f^i(E_i)$ is also in $\cs_{\alpha}$.
  Finally, since $f(t)\geq f^i(t)$ for every $ i$ and $ t$, then
  $f(E_i)\in\cs_{\alpha}\ (i=1,\dots,k)$ and since $f(1)\geq k$ we
  have $k\leq E_1<\dots<E_k$ so that $E=\cup_1^kE_i\in\cs_{\alpha+1}$
  as required.

  \textit{Case 3, $\alpha$ is a limit ordinal.} For the $\alpha$-game
  $\cn$ is bound to pick $l$ and then they play a bound
  $\alpha_l$-game (where $\alpha_n\nearrow\alpha$ is the sequence of
  ordinals increasing to $\alpha$ fixed in the definition of
  $\cs_{\alpha}$).  By assumption we may choose $f'$ for the
  $\alpha_l$-game such that if $E$ is the result of the
  $\alpha_l$-game where $\cs$ has chosen $E$ as small as possible then
  $f'(E)\in\cs_{\alpha_l}$. Let $f(t)=f'(t)+l$.  Now,
  $f'(E)\in\cs_{\alpha_l}$ for the $E$ we fixed initially, which
  implies $f(E)\in\cs_{\alpha_l}$, since $f(E)$ is a spreading of
  $f'(E)$.  Finally $f(1)\geq l$, hence $f(E)\in\cs_{\alpha}$ since
  $\{F\in\cs_{\alpha_l} : l\leq F \}\subseteq\cs_{\alpha}$.

  To generalize for $(\alpha_1,\dots,\alpha_r)$ we proceed as in
  Case~2, using bound $\alpha_i$-games $(i=1,\dots,r)$.
\end{proof}

\begin{cor} \label{C:F-in-S}
  If $\cn$ has a winning strategy for an
  $(\alpha_1,\dots,\alpha_r)$-Schreier game on a hereditary
  collection $\cf\subseteq[\N]^{<\omega}$, then there exists
  $N=(n_i)\in[\N]$ such that
  $\cf(N)\subseteq(\cs_{\alpha_1},\dots,\cs_{\alpha_r})$.
\end{cor}
\begin{proof}
  Suppose $\cn$ has a winning strategy for an
  $(\alpha_1,\dots,\alpha_r)$-Schreier game on
  $\cf\subseteq[\N]^{<\omega}$.  Let $\cn,\ \cs$ play the bound
  $(\alpha_1,\dots,\alpha_r)$-game where $\cn$ always chooses $l$ as
  small as possible so that $\cn$ will win.  Let $E\in\cf$, then we
  may decompose $E=\cup_1^pE_i$ according to this game as follows.  If
  the first set which $\cs$ chooses must have length greater than or
  equal to $l_1$, then let $E_1=\{e_1,\dots,e_{l_1}\}$; if $\cs$ has
  chosen $E_1<\dots<E_{q-1}$ and $\cs$ must pick the
  $q^{\mbox{\scriptsize th}}$ set to have length at least $l_q$, then
  let $E_q$ be the next $l_q$ elements of $E$ after $E_{q-1}$.  Since
  $\cn$ has a winning strategy, and $E\in\cf$, this process must
  exhaust $E$, but at that point let $\cs$ continue the game, always
  choosing sets as small as possible, and let $\bar E$ be the union of
  the sets obtained (including $E$).  Now, by
  Proposition~\ref{P:schreier}, there exists $N=(n_i)\in[\N]$ such
  that $n_F\in(\cs_{\alpha_1},\dots,\cs_{\alpha_r})$ for any set $F$
  resulting from such a game.  Thus $n_{\bar
  E}\in(\cs_{\alpha_1},\dots,\cs_{\alpha_r})$, and hence so is $n_E$
  since $E\subseteq\bar E$ and $\cs_{\alpha}$ is hereditary for each
  $\alpha$.  So finally,
  $\cf(N)\subseteq(\cs_{\alpha_1},\dots,\cs_{\alpha_r})$ as required.
\end{proof}

These results are sufficient to prove Theorem~\ref{bigthm}:

\begin{proof}[Proof of Theorem~\ref{bigthm}]
  Let $r\geq1$, let $0\leq\alpha_1\leq\dots\leq\alpha_r<\omega_1$, and
  let $\cf$ be a hereditary collection, $\cf\subseteq[\N]^{<\omega}$.
  Since $(\alpha_1,\dots,\alpha_r)$ has the Dichotomy property, it
  follows that there either exists $M\in[\bar N]$ such that
  $(\cs_{\alpha_1},\dots,\cs_{\alpha_r})(M)\subseteq\cf$, in which
  case the proof is complete, or there exists $M\in[\bar N]$ such that
  $\cn$ has a winning strategy for the
  $(\alpha_1,\dots,\alpha_r)$-Schreier game on $\cf[M]$.  Now $\cf[M]$
  is again hereditary and thus by Corollary~\ref{C:F-in-S} there
  exists $N\in[\N]$ such that
  $\cf[M](N)\subseteq(\cs_{\alpha_1},\dots,\cs_{\alpha_r})$ as
  required.
\end{proof}

Theorem~\ref{mainthm} follows from Theorem~\ref{bigthm} as an
immediate corollary.

\begin{rem}
It should be noted that Theorem~\ref{bigthm} is no longer true if we
do not first restrict $\cf$ to a subsequence of $\N$.  Indeed, we have
the following example:

\begin{expl}
  We construct a hereditary collection $\cf\subseteq[\N]^{<\omega}$
  such that for every $ M\in[\N]$ we have both $\cs_1(M)\nsubseteq\cf$
  and $\cf(M)\nsubseteq\cs_1$.
  
  Let $F_k=\{ 2^k+1,\dots,2^k+k \}$ for each $k\geq1$ and let
  \[\cf=\cup_{k=1}^{\infty}\{ \{1\}\cup E, E : E\subseteq F_k \}\ ,\]
  then $\cf$ is clearly hereditary.  Let $M\in[\N]$ and let $l=m_1$.
  Then $F=\{m_1\}\cup m_{F_l}\in\cf(M)$, but $|F|=l+1>\min F$ and
  hence $F\not\in\cs_1$.  Furthermore, $\cs_1(M)\nsubseteq\cf$ for
  suppose $E\in\cs_1(M)$ and $E\in\cf$ with $|E|>2$.  Let
  $E=\{e_1,\dots,e_p\}$ and find $k$ such that $E<F_k$.  Now let
  $E'=\{e_1,\dots,e_{p-1},m_{2^k+1}\}$, then $E'$ is still in
  $\cs_1(M)$ since this collection is spreading, but $E'\not\in\cf$
  because if $F,\ F'\in\cf$ then either $F\subseteq F',\ F'\subseteq
  F$ or $|F\cap F'|=0$ or $1$.  None of these is true for $E,E'$.
\end{expl}
\end{rem}

\section{Application}
\label{sec:app}

In this section we use Theorem~\ref{mainthm} to provide an alternative
proof of a result in a paper of Argyros, Mercourakis and
Tsarpalias~\cite{amt}.  We first state some definitions.

\begin{defn} \label{scbi} 
  Strong Cantor-Bendixson Index \cite{amt}
\end{defn}

Let $\cf$ be an adequate family (hereditary and closed) of finite
subsets of $\N$ as defined in Section~\ref{sec:prelims}.  For
$L\in[\N]$ we define the \emph{strong Cantor-Bendixson derivative} of
$\cf[L]$ for each ordinal $\alpha<\omega_1$ to be:
\[ \cf[L]^{(1)}=\{ A\in\cf[L] : A \mbox{ is a cluster point of } 
\cf[A\cup N]\ \mbox{ for each } N\in[L] \}\ .\] 
Thus for finite $A\subseteq L$
we have that $A\in\cf[L]^{(1)}$ if and only if $\{l\in L :
A\cup\{l\}\not\in\cf \}$ is finite.  If we have defined
$\cf[L]^{(\alpha)}$, the $\alpha^{\scriptstyle\textrm{th}}$ strong
Cantor-Bendixson derivative of $\cf[L]$, then we define the
$(\alpha+1)^{\scriptstyle\textrm{th}}$ derivative as:
\[ \cf[L]^{(\alpha+1)}=(\cf[L]^{(\alpha)})[L]^{(1)}\ .\]
If $\alpha$ is a limit ordinal and we have defined $\cf[L]^{(\beta)}$
for each $\beta<\alpha$, then we set
\[ \cf[L]^{(\alpha)} = \bigcap_{\beta<\alpha} \cf[L]^{(\beta)}\ .\]
The \emph{strong Cantor-Bendixson index} of $\cf[L]$ is defined to be
the smallest countable ordinal $\alpha<\omega_1$ such that
$\cf[L]^{(\alpha)}=\emptyset$.  We denote this index by $s(\cf[L])$.
For more detail concerning the strong Cantor-Bendixson derivative and
index please refer to~\cite{amt}.

\begin{rem} \label{amt-rems} 
  The following are stated in \cite{amt} or are simple consequences of
  their work:
  \begin{enumerate}
  \item The strong Cantor-Bendixson index must be a successor ordinal.
  \item For each $\alpha<\omega_1$ we have
    $s(\cs_{\alpha})=\omega^{\alpha}+1$ (\cite{amt}~Remark~2.2.5).
  \item If $\cf\subseteq[\N]^{<\omega}$ is spreading, then
    $s(\cf[L])=s(\cf)$ for every $L\in[\N]$.
  \item If $s(\cf[L])>\alpha$, then $s(\cf[M])>\alpha$ for every
    $M\in[L]$ (\cite{amt}~Proposition~2.2.3).
  \end{enumerate}
\end{rem}

We prove the following result from~\cite{amt} (Theorem~2.2.6):

\begin{thm} \label{amt}
  \emph{(\cite{amt})} Let $\cf\subseteq[\N]^{<\omega}$ be an adequate
  family.  If there exists $L\in[\N]$ such that
  $s(\cf[L])>\omega^{\alpha}$, then there exists $M\in[L]$ such that
  $\cs_{\alpha}(M)\subseteq\cf[M]$.
\end{thm}

\begin{proof} 
  Let $\cf\subseteq[\N]^{<\omega}$ be an adequate family and let
  $L\in[\N]$ satisfy $s(\cf[L])>\omega^{\alpha}$.  Suppose first that
  in fact $s(\cf[L])>\omega^{\alpha}+1$.  Now, by
  Theorem~\ref{mainthm}, either there exists $M\in[L]$ such that
  $\cs_{\alpha}(M)\subseteq\cf[M]$ as required, or else there exist
  $M\in[L],\ N\in[\N]$ such that $\cf[M](N)\subseteq\cs_{\alpha}$.  We
  can easily see that the index of $\cf[M]$ is the same as the index
  of $(\cf[M](N))[N_M]$ where $N_M=(n_m)_{m\in M}$.  Indeed, if
  $(A_i)$ is a sequence in $\cf[M]^{(\beta)}$ converging to
  $A\in\cf[M]^{(\beta+1)}$, then $(n_{A_i})$ is a sequence in
  $(\cf[M](N))[N_M]^{(\beta)}$ converging to
  $n_A\in(\cf[M](N))[N_M]^{(\beta+1)}$ and vice versa.  Thus, if
  $\cf[M](N)\subseteq\cs_{\alpha}$, then \[ s(\cf[M]) =
  s((\cf[M](N))[N_M]) \leq s(\cs_{\alpha}[N_M]) = \omega^{\alpha}+1\ .
  \] However, by Remark~\ref{amt-rems}~(iv),
  $s(\cf[M])>\omega^{\alpha}+1$, a contradiction.  Thus the second
  case, above, cannot happen.

  To finish the proof we assume that $s(\cf[L])=\omega^{\alpha}+1$ and
  define
  \[ \bar{\cf}=\{ \{n\}\cup F : F\in\cf,\ n<F \}\cup\cf\ . \]
  If $A\in\cf[L]^{(\beta)}\setminus\{\emptyset\}$ and $l\in L$ with
  $l<A$, then $\{l\}\cup A\in\bar{\cf}[L]^{(\beta)}$.  Indeed, suppose
  this is true for some ordinal $\beta<\alpha$ and let
  $A\in\cf[L]^{(\beta+1)}\setminus\{\emptyset\}$ and $l\in L$ with
  $l<A$, then there exists a sequence $(A_i)\subseteq\cf[L]^{(\beta)}$
  converging to $A$.  Now, $(\{l\}\cup A_i)$ is a sequence in
  $\bar{\cf}[L]^{(\beta)}$ converging to $\{l\}\cup A$, hence
  $A\in\bar{\cf}[L]^{(\beta+1)}$ as required.  The limit ordinal case
  is clear.  Since $\cf[L]^{(\omega^{\alpha})}\ne\emptyset$, it
  follows that $\cf[L]^{(\beta)}$ is infinite for each
  $\beta<\omega^{\alpha}$, so that $\{l\}\in\bar{\cf}[L]^{(\beta)}$
  for every $ l\in L$, and each $\beta<\omega^{\alpha}$.  Thus
  $\{l\}\in\bar{\cf}[L]^{(\omega^{\alpha})}$ for every $ l\in L$ and
  hence $\emptyset\in\bar{\cf}[L]^{(\omega^{\alpha}+1)}$, so that
  $s(\bar{\cf}[L])>\omega^{\alpha}+1$.  Finally, we apply the previous
  case to $\bar{\cf}[L]$ to obtain $M=(m_i)\in[L]$ with
  $\cs_{\alpha}(M)\subseteq\bar{\cf}[M]$, then setting
  $M'=(m_i)_{i>2}$ we have $\cs_{\alpha}(M')\subseteq\cf[M']$ as
  required.
\end{proof}

\end{document}